\documentclass[twoside,11pt,leqno]{article}
\usepackage{amsfonts}

 \def\X{{\cal X}}  \def\H{{\cal H}} 
\def\D{{\texttt{D}}}

\def\B{B({\cal H})} 
\def\asc{ \textrm{asc}} \def\dsc{ \textrm{dsc}}

\newtheorem{df}{Definition}[section]
\newtheorem{thm}[df]{Theorem} \newtheorem{pro}[df]{Proposition}
 
\newtheorem{rema}[df] {Remark} \newtheorem{lem}[df] {Lemma}
\def\sfstp{{\hskip-1em}{\bf.}{\hskip1em}}

\def\subject#1{\renewcommand{\thefootnote}{}\footnote
{AMS(MOS) subject classification (2010). Primary: {#1}}}

\def\keywords#1{\renewcommand{\thefootnote}{}\footnote
{Keywords: {#1}}}

\def\enddemo{\qed \endtrivlist} \expandafter\let\csname
enddemo*\endcsname=\enddemo

\def\qedsymbol{\ifmmode\bgroup\else$\bgroup\aftergroup$\fi
\vcenter{\hrule\hbox{\vrule
height.5em\kern.5em\vrule}\hrule}\egroup}
\def\qed{\ifmmode\else\unskip\nobreak\fi\quad\qedsymbol}

\title{\bf  Drazin invertible $(m, P)$-expansive operators}
\author{\normalsize B.P. Duggal, I.H. Kim}
\date{}

\begin{document}

\maketitle \thispagestyle{empty} \vskip-16pt

\subject{Primary47A63 Secondary47A11.} 
\footnote{The second named author was supported by Basic Science Research Program through the National Research Foundation of Korea (NRF)
 funded by the Ministry of Education (NRF-2019R1F1A1057574)}
\begin{abstract} A Hilbert space operator $T\in\B$ is $(m,P)$-expansive, for some positive integer $m$ and operator $P\in\B$, if $\sum_{j=0}^m{(-1)^j\left(\begin{array}{clcr}m\\j\end{array}\right)T^{*j}PT^j}\leq 0$. No Drazin invertible operator $T$
can be $(m,I)$-expansive, and if $T$ is $(m,P)$-expansive for some positive operator $P$, then necessarily $P$ has a decomposition $P=P_{11}\oplus 0$. If $T$ is $(m,|T^n|^2)$-expansive for some positive integer $n$, then $T^n$ has a decomposition $T^n=\left(\begin{array}{clcr}U_1P_1 & X\\0 & 0\end{array}\right)$; if also $\left(\begin{array}{clcr}I_1 & X\\X^* & X^*X\end{array}\right)\geq I$, then $\left(\begin{array}{clcr}P_1U_1 & P_1X\\0 & 0\end{array}\right)$ is $(m,I)$-expansive and $\left(\begin{array}{clcr}P^{\frac{1}{2}}_1U_1P^{\frac{1}{2}}_1 & P_1^{\frac{1}{2}}X\\0 & 0\end{array}\right)$ is $(m,I)$-expansive in an equivalent norm on $\H$.

\end{abstract}


\section {\sfstp Introduction} Let  $\B$ denote the algebra of operators, with identity $I$,  on an infinite dimensional separable complex  Hilbert space $\H$  into itself.  For $A, B\in\B$, let $L_A$ and $R_B$ $\in B(\B)$ denote, respectively, the operators
$$
L_A(X)=AX, \ R_B(X)=XB
$$
of left multiplication by $A$ and right multiplication by $B$. For an operator $X\in\B$ and an integer $m\geq 1$, $T\in \B$ is an $(m,X)$-expansive operator, denoted $T\in (m,X)$-expansive, if $$\triangle^m_{T^*,T}(X)=(I-L_{T^*}R_T)^m(X)=\sum_{j=0}^m{(-1)^j\left(\begin{array}{clcr}m\\j\end{array}\right)T^{*j}XT^j}\leq 0.$$ The concept of $(m,X)$-expansive operators is an extension of the concept of $(m,I)$-expansive, equivalently $m$-expansive operators, and there is a large body of information available in
extant literature on $m$-expansive operators and their  various extensions (see \cite{{Ah}, {At}, {EJL}, {G}, {JKKL}, {SA}} for further references). A particularly rich theory is available for $(m,I)$-expansive operators $T$ for which $\triangle^m_{T^*,T}(I)=0$: Such  operators have been called $m$-isometric operators -- see \cite
{{AS1}, {FB},  {BMN},  {Fb}, {DM}} for further references.

\vskip4pt Let $[S,T]=ST-TS$ denote the commutator of operators $S,T\in\B$. $T\in\B$ is Drazin invertible, with Drazin inverse $T_d$, if 
$$
 [T_d,T]=0, \ T_d^2T=T_d, \  T^{p+1}T_d=T^p
 $$ 
for some integer $p\geq 1$. The least integer $p$ is then said to be the Drazin index of $T$. We start this note by proving
that Drazin invertible operators $T\in\B$ can not be $m$-expansive, and for this we give a short, simple  proof (using little more than the elementary properties of left/right multiplication oprerators and some basic algebra) of a result from \cite{EJL} which says that $T\in m$-expansive implies $T^n\in m$-expansive for all integers $n\geq 1$. There may however exist Drazin invertible $T$ such that $T\in (m,P)$-expansive for some $P\geq 0$. For such $T$ and $P$ it is seen that $P$ has  a decomposition $P=P_{11}\oplus 0$ and $(-1)^m\triangle_{T_d^*,T_d}(P)\leq 0$.  In the particular case in which $m=2$, this implies that $T$ is $P$-isometric (i.e., $T^*PT=P$). For $T\in (m,P)$-expansive such that $P=T^{*n}T^n$ for some positive integer $n$, $T^n$ has an upper triangular matrix representation $T^n=\left(\begin{array}{clcr} T^n_1 & X\\ 0 & 0\end{array}\right)$, where $T_1$ is $m$-expansive. Letting $T^n_1$ have the polar decompositon
$T^n_1=U_1P_1$ and assuming  $\left(\begin{array}{clcr}I_1 & X\\X^* & X^*X\end{array}\right)\geq I$, we prove that the operator $ \left(\begin{array}{clcr}P_1U_1 & P_1X \\0 & 0\end{array}\right)$ is $m$-expansive and that there exists an equivalent norm on $\H$ in which the operator $\left(\begin{array}{clcr} P_1^{\frac{1}{2}}U_1P_1^{\frac{1}{2}}& P_1^{\frac{1}{2}}X\\0 & 0\end{array}\right)$ is $m$-expansive.

\section {\sfstp Drazin invertible $m$-expansive operators}  If an operator $T\in\B$ is Drazin invertible, then $\asc(T)=\dsc(T)\leq p$ for some positive integer $p$ and there exists a direct sum decomposition
$$
\H=T^p(\H)\oplus T^{-p}(0)=\H_1\oplus\H_2,  \hspace{2mm}T^{-p}(0)=T^{{p}^{-1}}(0)
$$
of $\H$ such that
$$
T=T|_{\H_1}\oplus T|_{\H_2}=T_1\oplus T_2\in B(\H_1\oplus\H_2),
$$
where $T_1$ is invertible and $T_2$ is $p$-nilpotent. The Drazin inverse $T_d$ of $T$ is then the operator
$$
T_d=T^{-1}_1\oplus 0\in B(\H_1\oplus\H_2)
$$
\cite[Theorem 2.2.3]{DR}. No Drazin invertible operator can be $m$-expansive (for any positive integer $m$). To see this we start with a proof of the following well known result from \cite{EJL}, Lemmas 2.1, 2.2 and Theorem 2.3: Our proof here is much shorter and different from that of  \cite{EJL}.

\begin{lem}\label{lem0} $T\in\B$ $m$-expansive implies $T^n\in m$-expansive for all integers $n\geq 1$.
\end{lem}

\begin{demo}$T$ is $m$-expansive if and only if $\triangle^m_{T^*,T}(I)\leq 0$; hence
\begin{eqnarray*} \triangle^m_{T^{*n},T^n}(I ) &=& (I-L_{T^{*n}}R_{T^n})^m(I)=(I-L^n_{T^*}R^n_T)^m(I)\\
&=&\{L^{n-1}_{T^*}\triangle_{T^*,T}R^{n-1}_T+L^{n-2}_{T^*}\triangle_{T^*,T}R^{n-2}_T+ \cdots + \\
&+& L_{T^*}\triangle_{T^*,T}R_T+\triangle_{T^*,T}\}^m(I)\\ &=& \{L^{n-1}_{T^*}R^{n-1}_T+L^{n-2}_{T^*}R^{n-2}_T+ \cdots + L_{T^*}R_T+ I \}^m\triangle_{T^*,T}^m(I)\\
&=&  \sum_{j=0}^{m(n-1)}{\alpha_{m(n-1)-j}L_{T^*}^{m(n-1)-j}\triangle^m_{T^*,T}(I)R_T^{m(n-1)-j}}\\
&\leq& 0,
\end{eqnarray*}
since $\alpha_{m(n-1)-j}\geq 0$ for all $j$.
\end{demo}

It is immediate from the argument above that $T\in (m,X)$-expansive for some operator $X\in\B$ implies $T^n\in (m,X)$-expansive for all integers $n\geq 1$.

\begin{thm}\label{thm0} There exist no Drazin invertible $m$-expansive operators $T\in\B$.
\end{thm}

\begin{demo} Assume that there exists a Drazin invertible $m$-expansive operator $T$. Let $p$ be the Drazin index of $T$. Then $T^p\in m$-expansive (by Lemma \ref{lem0}), and we have
\begin{eqnarray*} & & \sum_{j=0}^m{(-1)^j\left(\begin{array}{clcr}m\\j\end{array}\right) T^{*pj}T^{pj}}\leq 0\\
&\Longleftrightarrow& \sum_{j=0}^m{(-1)^j\left(\begin{array}{clcr}m\\j\end{array}\right) (T_1^*\oplus T_2)^{pj}(T_1\oplus T_2)^{pj}}\leq 0\\
&\Longleftrightarrow& \left(\sum_{j=0}^m{(-1)^j\left(\begin{array}{clcr}m\\j\end{array}\right) (T_1^*)^{pj}(T_1)^{pj}}\right) \oplus \left(\sum_{j=0}^m{(-1)^j\left(\begin{array}{clcr}m\\j\end{array}\right) (T_2^*)^{pj}(T_2)^{pj}}\right)\leq 0\\
&\Longleftrightarrow&  \left( \sum_{j=0}^m{(-1)^j\left(\begin{array}{clcr}m\\j\end{array}\right) T_1^{*pj}T_1^{pj}}\right)  \oplus  I_2\leq 0,
\end{eqnarray*}
where $I_2$ is the identity of $B(\H_2)$. The final inequality obviously being false, the theorem is proved.
\end{demo}

There may however exist operators $X\in\B$ such that $\triangle^m_{T^*,T}(X)\leq 0$, i.e. $T\in (m,X)$-expansive, for Drazin invertible $T=T_1\oplus T_2\in B(\H_1\oplus\H_2)$. For example, if $\triangle^m_{T^*_1,T_1}(I_1)\leq 0$, $I_1$ the identity of $B(\H_1)$,
then $\triangle^m_{T^*,T}(X)\leq 0$ for $X=I_1\oplus 0\in B(\H_1\oplus\H_2)$.The following theorem asserts that if $T\in (m,P)$-expansive for some positive operator $P\in\B$, then necessarily $P=P_{11}\oplus 0\in B(\H_1\oplus\H_2)$.  Let
$$
{\tilde\triangle}^m_{T^*,T}(X)=(L_{T^*}R_T-I)^m(X)=(-1)^m\triangle^m_{T^*,T}(X).
$$

\begin{thm}\label{thm00}  A Drazin invertible operator $T\in B(\H_1\oplus\H_2)$ with Drazin inverse $T_d$ is $(m,P)$-expansive for some positive operator $P\in\B$ if and only if $P=P_{11}\oplus 0\in B(\H_1\oplus\H_2)$ and ${\tilde\triangle}^m_{T_d^*,T_d}(P)\leq 0$. 
\end{thm}

\begin{demo} Letting $P\in B(\H_1\oplus\H_2)$ have the matrix representation $$P=[P_{ij}]_{i,j=1,2},$$ the hypothesis  $T=T_1\oplus T_2$, $T_1$ invertible and $T^p_2=0$ for some positive integer $p$, is $(m,P)$-expansive implies
  \begin{eqnarray}  
  \triangle^m_{T^*,T}(P)=\left[ \sum_{k=0}^m{(-1)^k\left(\begin{array}{clcr}m\\k\end{array}\right) T_i^{*k}P_{ij}T_j^{k}} \right]_{i,j=1,2}\leq 0. 
  \end{eqnarray}
We prove that $P_{22}=0$. Inequality $(1)$ implies
\begin{eqnarray} 
\sum_{k=0}^m{(-1)^k\left(\begin{array}{clcr}m\\k\end{array}\right) T_2^{*k}P_{22}T_2^{k}}\leq 0,
\end{eqnarray}
and hence upon multiplying by $T_2^{*(p-1)}$ on the left and by $T_2^{p-1}$ on the right (in inequality $(2)$)
$$
\sum_{k=0}^m{(-1)^k\left(\begin{array}{clcr}m\\k\end{array}\right) T_2^{*(k+p-1)}P_{22}T_2^{k+p-1}}\leq 0.
$$
Since $T_2^p=0$, we have 
$$
T_2^{*(p-1)}P_{22}T_2^{p-1}\leq 0.
$$
Considering now  
$$ 
T_2^{*(p-2)} \left( \sum_{k=0}^m{(-1)^k\left(\begin{array}{clcr}m\\k\end{array}\right) T_2^{*(k)}P_{22}T_2^{k}} \right) T_2^{p-2},
$$
we have 
$$
T_2^{*(p-2)}P_{22}T^{p-2}-mT_2^{*(p-1)}P_{22}T_2^{p-1}\leq 0.
$$ 
Consequently, 
$$T_2^{*(p-2)}P_{22}T^{p-2}\leq mT_2^{*(p-1)}P_{22}T_2^{p-1}.
$$
Repeating this argument, multiplying this time by $T_2^{*(p-3)}$ on the left and by $T_2^{p-3}$ on the right, we have
\begin{eqnarray*}
 & & \sum_{k=0}^m{(-1)^k\left(\begin{array}{clcr}m\\k\end{array}\right) T_2^{*(k+p-3)}P_{22}T_2^{k+p-3}}\leq 0\\
&\Longleftrightarrow& T_2^{*(p-3)}P_{22}T_2^{p-3}-mT_2^{*(p-2)}P_{22}T_2^{p-2}+{\frac{m(m-1)}{2}}T_2^{*(p-1)}P_{22}T_2^{p-1}\leq 0\\
&\Longrightarrow&  T_2^{*(p-3)}P_{22}T_2^{p-3}\leq mT_2^{*(p-2)}P_{22}T_2^{p-2}-{\frac{m(m-1)}{2}}T_2^{*(p-1)}P_{22}T_2^{p-1}\\
&\Longrightarrow& T_2^{*(p-3)}P_{22}T_2^{p-3}\leq \left(m^2-{\frac{m(m-1)}{2}}\right)T_2^{*(p-1)}P_{22}T_2^{p-1}\\
&\Longrightarrow& {\frac{m(m+1)}{2}}T_2^{*(p-1)}P_{22}T_2^{p-1}\leq 0.
\end{eqnarray*}
An induction argument now shows that indeed
$$
T_2^{*(p-r)}P_{22}T_2^{p-r}\leq {\frac{m(m+1) \cdots (m+r-2)}{(r-1)!}}T_2^{*(p-1)}P_{22}T_2^{p-1}\leq 0.
$$
Hence 
$$
P_{22}=(I_2-L_{T^*_2}R_{T_2})^m(P_{22})\leq 0.
$$
Recall now from \cite{An}, Theorem I.1, that an operator $P=[P_{ij}]_{i,j=1,2}\geq 0$ if and only if $P_{11}, P_{22}\geq 0$ and there exists a contraction $C$ such that $P_{12}^*=P_{21}=P_{22}^{\frac{1}{2}}CP_{11}^{\frac{1}{2}}$. Since $P_{22}\geq 0$
forces $P_{22}=0$, and then $P_{21}=0$ for every contraction $C\in B(\H_1,\H_2)$, we conclude that $P=P_{11}\oplus 0\in B(\H_1,\H_2)$. Consequently,
$$
T\in (m,P)-{\rm expansive}\Longleftrightarrow T_1\in (m,P_{11})-{\rm expansive}.
$$
 The operator $T_1$ being invertible, 
 \begin{eqnarray*}
  \triangle^m_{T_1^*,T_1}(P_{11})\leq 0 &\Longleftrightarrow& T^{*m}_1\tilde\triangle^m_{T^{-*1}_1T_1^{-1}}(P_{11})T^m_1\leq 0 \\
&\Longleftrightarrow& \tilde\triangle^m_{T^{-*1}_1T_1^{-1}}(P_{11})\leq 0\\
&\Longleftrightarrow& \tilde\triangle^m_{T_d^*,T_d}(P)\leq 0, \  P=P_{11}\oplus 0.
\end{eqnarray*} 
This completes the proof.
\end{demo}

Observe that $\tilde\triangle^m_{T^*,T}(P)=\triangle^m_{T^*,T}(P)$ in the case in which $m$ is even. The case $m=2$ is of some interrest.  Remark here that the introduction of the weight $P$, i.e. the non-negative operator $P\in\B$, results in an
additional semi-inner product ${\langle x,y \rangle}_P=\langle Px,y\rangle$ on $\H$ for every $x,y\in\H$. The additional structure gives rise to an adjoint operator, not defined for every $A\in\B$ unless $P$ is invertible. Let $||.||_P$ denote the semi-norm induced by $\langle , \rangle_P$, i.e., $||P^{\frac{1}{2}}x||^2=||x||_P^2=\langle x,x\rangle_P$.
Then $||.||_P$ is a norm on $\H$ if and only if $P$ is injective. If $P$ is invertible, then $||P^{-\frac{1}{2}}||^{-2}||x||^2\leq ||P^{\frac{1}{2}}x||^2\leq ||P^{\frac{1}{2}}||^2||x||^2$, and $||.||_P$ defines an equivalent norm on $\H$. (See \cite{ACG} for further information.)
We say in the following that an operator $A\in\B$ is $P$-isometric for some positive operator $P\in\B$,
if $||Ax||^2_P=\langle A^*PAx,x\rangle=||P^{\frac{1}{2}}Ax||^2=||x||^2_P$ for every $x\in\H$  \cite[Definition 3.1]{ACG}.

\begin{pro}\label{pro01} If $T\in\B$ is a Drazin invertible $(2,P)$-expansive operator for some positive operator $P\in\B$, then $T$ is $P$-isometric.
\end{pro}

\begin{demo} The hypothesis
\begin{eqnarray*} & & T\in (2,P)-{\rm expansive} \Longleftrightarrow \triangle^2_{T^*,T}(P)\leq 0\\
&\Longleftrightarrow& T^{*2}PT^2\leq 2T^*PT-P\Longleftrightarrow T^{*2}PT^2\leq 2\tilde\triangle_{T^*,T}(P)+P\\
&\Longrightarrow& T^{*3}PT^3\leq 2T^*\tilde\triangle_{T^*,T}(P)+T^*PT\leq 3\tilde\triangle_{T^*,T}(P)+P,
\end{eqnarray*}
and (hence) by an induction argument that $$T^{*n}PT^n\leq n\tilde\triangle_{T^*,T}(P)+P\Longleftrightarrow \tilde\triangle_{T^*,T}(P)\geq {\frac{1}{n}}\{T^{*n}PT^n-P\}.$$
Thus
\begin{eqnarray*}  \langle\tilde\triangle_{T^*,T}(P)x,x\rangle \hspace{3mm}&\geq & {\lim\sup}_{n} {\frac{1}{n}}\langle (T^{*n}PT^n-P)x,x\rangle \\
&\Longrightarrow& \tilde\triangle_{T^*,T}(P)\geq 0\\ &\Longleftrightarrow& T^*PT\geq P.
\end{eqnarray*}
The operator $T$ being Drazin invertible, we also have  $$\tilde\triangle^2_{T^*_d,T_d}(P)=\triangle^2_{T^*_d,T_d}(P)\leq 0,$$
and this upon arguing as above implies 
\begin{eqnarray*} T^*_dPT_d\geq P &\Longleftrightarrow& {T_1^*}^{-1}P_{11}T^{-1}_1\geq P_{11}, \ P=P_{11}\oplus 0\\
&\Longleftrightarrow& P_{11}\geq T^*_1P_{11}T_1, \ P=P_{11}\oplus 0\\
&\Longleftrightarrow& P\geq T^*PT, \ P=P_{11}\oplus 0. 
\end{eqnarray*}
Consequently, $T^*PT=P$, which is equivalent to $T$ is a $P$-isometry.\end{demo}

\begin{rema}\label{rema01} It is immediate from the argument of the proof of Proposition \ref{pro01} that if $T$ is an invertible $(2,P)$-expansive operator, then $T$ is a $P$-isometry.
If $T$ is Drazin invertible and $(2,|T|^2)$-expansive, then $T_1$ is unitary and $T$ is the direct sum of a unitary and a nilpotent.
\end{rema}

Hiding within the argument of the proof of Proposition \ref{pro01} is a more general result. Recall that $T\in\B$ is $m$-contractive, equivalently $(m,I)$-contractive, if $\triangle_{T^*,T}^m(I)\geq 0$. Defining $(m,P)$-contractive, $P\in\B$, by requiring $\triangle^m_{T^*,T}(P)\geq 0$, it
is clear that operators $T\in \{(m,P)-{\rm expansive}\} \wedge  \{(m,P)-{\rm contractive}\}$ are $(m,P)$-isometric. Let $m\geq 2$, and consider an $(m,P)$-expansive operator $T$, $P\in \B$ some positive operator, such that $T$ is (also) $(m-2,P)$-contractive.
(The $(m-2,P)$-contractive hypothesis on $T$ is vacuous in the case in which $m=2$.) Define the positive operator
$Q\in\B$ by $Q=\triangle^{m-2}_{T^*,T}(P)$. Then $\triangle^2_{T^*,T}(Q)\leq 0$, and it follows from the argument of the proof of Proposition \ref{pro01} that 
$$
T^*QT\geq Q \ (\Longleftrightarrow \ \triangle^{m-1}_{T^*,T}(P)\leq 0).
$$ 
Assume now additionally that $T$ is Drazin invertible.
Then, see the proof of Proposition \ref{pro01}, $T^*QT\leq Q$, and hence 
$$
\triangle^{m-1}_{T^*,T}(P)=0.
$$
It is easily verified that $T\in (n,P)$-isometric for some positive integer $n$ implies $T\in (t,P)$-isometric for all integers $t\geq n$. (Indeed, $\triangle^t_{T^*,T}(P)=\triangle^{t-n}_{T^*,T}(\triangle^{n}_{T^*,T}(P))$.)  Hence $\triangle^m_{T^*,T}(P)=0$. We have proved:

\begin{thm}\label{thm01} If a Drazin invertible operator $T\in\B$ is both $(m,P)$-expansive and $(m-2,P)$-contractive for some positive operator $P\in\B$ and positive integer $m\geq 2$, then $T\in (m-1,P)$-isometric.
\end{thm}

\section {\sfstp  $(m,P)$-expansive operators: $P=T^{*n}T^n$}  If $T\in (m,P)$-expansive for a positive operator $P$ such that $0$ is not in the approximate point spectrum $\sigma_a(P)$ of $P$, and $\{x_n\}\subset\H$ is a sequence of unit vectors such that
$$
||(T-\lambda)x_n||\longrightarrow 0 \ \mbox {as} \ n\longrightarrow\infty,
$$
then
\begin{eqnarray*} \sum_{j=0}^m{(-1)^j\left(\begin{array}{clcr}m\\j\end{array}\right)\langle T^{*j}PT^jx_n,x_n\rangle } & = &  \sum_{j=0}^m{(-1)^j\left(\begin{array}{clcr}m\\j\end{array}\right)|\lambda|^{2j}||P^{\frac{1}{2}}x_n||^2}\\
&=& (1-|\lambda|^2)^m||P^{\frac{1}{2}}x_n||^2\\ &\leq& 0 
\end{eqnarray*}
for all $n$. This, since $\lambda=0$ forces $\lim_{n\rightarrow\infty}||P^{\frac{1}{2}}x_n||=0$, implies $0\notin\sigma_a(T)$; again, if $\lambda\neq 0$ and $m$ is even then $|\lambda|=1$, and if $\lambda\neq 0$ and $n$ is odd then $|\lambda|\geq 1$.
Consequently, since $||T||\geq r(T)$ ($=$ the spectral radius $\lim_{n\rightarrow\infty}||T^n||^{\frac{1}{n}}$) for every operator $T$, $||T||\geq 1$.

\begin{rema}\label{rema00} The conclusion above  that $0\notin\sigma_a(T)$ (and $\sigma_a(T)$ is contained in a circle centered at the origin of radius greater than or equal to one) for operators $T\in (m,I)$-expansive ensures that $T$ can not be Drazin invertible. 
\end{rema}

\vskip4pt

Choose $P$ to be the positive operator $P=T^{*n}T^n$ for some integer $n\geq 1$, and assume $T\in (m,|T^n|^2)$-expansive. Define operators $T_1\in B(\overline{T^n(\H)})$ and $T_2\in B(T^{-*n}(0))$ by
$$ 
T_1=T|_{\overline{T^n(\H)}} \  \mbox{and} \ T_2=T|_{T^{-*n}(0)};
$$
then $T$ has an upper triangular matrix representation
$$
T=\left( \begin{array}{clcr} T_1 & *\\0 & T_2\end{array}\right) \in B(\overline{T^n(\H)}\oplus T^{-*n}(0)),
$$
where 
$$
T_1\in (m,|T_1^n|^2)-{\rm expansive}, \ T^n_2=0.
$$
The operator $T^n_1$ has a dense range; hence 
$$
T_1\in (m,|T^n_1|^2)-{\mbox{expansive}}\Longleftrightarrow T_1\in m-{\rm expansive}.
$$
Let $T^n_1$ have the polar decomposition $$T_1^n=U_1P_1;$$ then (the fact that $0\notin\sigma_a(T_1)$ implies) $U_1$ is an isometry and $P_1\geq 0$ is invertible. Define operators $A_1$ and $A_2$ by
$$
A_1=P^{\frac{1}{2}}U_1P^{\frac{1}{2}} \ {\rm and} \ A_2=P_1U_1.
$$
(Such operators $A_1$, respectively $A_2$, have been called the Aluthge transform \cite{{Al}, {FJKP}}, respectively the Duggal transform \cite{FJKP}, of $T^n_1=U_1P_1$.)  The operator $T_2$ being $n$-nilpotent, $T^n$ has an upper triangular representation
$$
T^n=\left(\begin{array}{clcr}T^n_1 & X\\0 & 0\end{array}\right).
$$
 Define operators $A, B, C, D\in  B(\overline{T^n(\H)}\oplus T^{-*n}(0))$ by
\begin{eqnarray*}
 & & A=\left(\begin{array}{clcr}A_1 & P_1^{\frac{1}{2}}X\\0 & 0\end{array}\right), \  B=\left(\begin{array}{clcr}A_2 & P_1X \\0 & 0\end{array}\right)\\
& & C=\left(\begin{array}{clcr}P_1 & P_1^{\frac{1}{2}}U^*_1X\\X^*U_1P_1^{\frac{1}{2}} & X^*X\end{array}\right) \ {\rm and} \  D=\left(\begin{array}{clcr}I_1 & U^*_1X\\X^*U_1 & X^*X\end{array}\right).
\end{eqnarray*}
If we now let 
$$
Q_1=P^{\frac{1}{2}}_1\oplus I_2, \ I_2  \ {\rm the  \ identity \ of} \ B(T^{-*n}(0)),
$$ then 
$$
A=Q^{-1}_1BQ_1, \ C=Q_1DQ_1  \ {\rm and} \ D=\left(\begin{array}{clcr}U^*_1 & 0\\X^* & 0\end{array}\right)\left(\begin{array}{clcr}U_1 & X\\0 & 0 \end{array}\right) \geq 0.
$$
The following theorem gives conditions ensuring $T\in (m,|T^n|^2)$-expansive implies $A$ and/or $B$ is $(m,P)$-expansive. As we shall see here, the operator $X$ in the upper triangular represetation of $T^n$ plays an important role. Let the operators $A, B, C, D, X, U_1, P_1$ and $Q_1$
be defined as in the above.

\begin{thm}\label{thm10} If $T\in (m,|T^n|^2)$-expansive, then $A\in (m,C)$-expansive and $B\in (m,D)$-expansive; if also $\left(\begin{array}{clcr}I_1 & X\\X^* & X^*X\end{array}\right)\geq I$, then $B\in m$-expansive and there exists an equivalent norm $||.||_Q$ on $\H$ such that
$A\in m$-expansive in this norm (equivalently, $A\in (m,Q)$-expansive). 
\end{thm}

\begin{demo} The hypothesis $T\in (m,|T^n|^2)$-expansive implies 
$$
T^{*n}\left[\sum_{j=0}^m{(-1)^j\left(\begin{array}{clcr}m\\j\end{array}\right)T^{*tj}T^{tj}}\right]T^n\leq 0
$$ 
for all integers $t\geq 1$. Choosing $t=n$ and letting $T^n=\left(\begin{array}{clcr}U_1P_1 & X\\0 & 0\end{array}\right)$
as above, this implies
\begin{eqnarray*}
 & &   \left(\begin{array}{clcr}P_1U^*_1 & 0\\X^* & 0\end{array}\right) \left[ \sum_{j=0}^m{(-1)^{j}\left(\begin{array}{clcr}m\\j\end{array}\right)\left(\begin{array}{clcr}P_1U_1^* & 0\\X^* & 0\end{array}\right)^j \left(\begin{array}{clcr}U_1P_1 & X\\0 & 0\end{array}\right)^{j}}\right]\times\\
&\times&\left(\begin{array}{clcr}U_1P_1 & X\\0 & 0\end{array}\right)\leq 0\\
 &\Longleftrightarrow& \left(\begin{array}{clcr} P_1 &0\\0 & I_2\end{array}\right)\left[ \sum_{j=0}^m{(-1)^{j}\left(\begin{array}{clcr}m\\j\end{array}\right)\left(\begin{array}{clcr}U^*_1P_1 & 0\\X^*P_1 & 0\end{array}\right)^j}
 \left(\begin{array}{clcr}U^*_1 &  0\\X^* & 0\end{array}\right)\left(\begin{array}{clcr}U_1 & X\\0 & 0\end{array}\right)\right.\times\\
 &\times&\left.\left(\begin{array}{clcr}P_1U_1 & P_1X\\0 & 0\end{array}\right)^j\right]\left(\begin{array}{clcr}P_1 & 0\\0 & I_2\end{array}\right) \leq 0\\
&\Longleftrightarrow&  \left(\begin{array}{clcr} P_1^{\frac{1}{2}} &0\\0 & I_2\end{array}\right) \sum_{j=0}^m(-1)^{j}\left(\begin{array}{clcr}m\\j\end{array}\right)\left(\begin{array}{clcr}P_1^{\frac{1}{2}}U^*_1P_1^{\frac{1}{2}} & 0\\X^*P_1^{\frac{1}{2}} & 0\end{array}\right)^j \left(\begin{array}{clcr}P_1^{\frac{1}{2}}U^*_1 & 0\\X^* & 0\end{array}\right)\times\\ 
&\times &\left(\begin{array}{clcr}U_1P_1^{\frac{1}{2}} & X\\0 & 0\end{array}\right) {\left(\begin{array}{clcr}P_1^{\frac{1}{2}}U_1P_1^{\frac{1}{2}} & P_1^{\frac{1}{2}}X\\0 & 0\end{array}\right)^j } \left(\begin{array}{clcr} P_1^{\frac{1}{2}} & 0\\0 & I_2\end{array}\right) \leq 0
\end{eqnarray*}
The operator $P_1^{\frac{1}{2}}+I_2=Q_1$ being invertible, we have
$$
T\in(m,|T^n|^2)-{\rm expansive}\Longrightarrow \sum_{j=0}^m{(-1)^{j}\left(\begin{array}{clcr}m\\j\end{array}\right) B^{*j}DB^j}\leq 0
$$ 
and 
$$ 
T\in(m, |T^n|^2)-{\rm expansive}\Longrightarrow \sum_{j=0}^m{(-1)^{j}\left(\begin{array}{clcr}m\\j\end{array}\right) A^{*j}CA^j}\leq 0.
$$
Observe that
\begin{eqnarray*}
 A^*CA &=& \left(\begin{array}{clcr}P_1^{\frac{1}{2}}U^*_1P_1^{\frac{1}{2}} & 0\\X^*P_1^{\frac{1}{2}} & 0\end{array}\right) \left(\begin{array}{clcr}P_1 & P_1^{\frac{1}{2}}U^*_1X\\X^*U_1P_1^{\frac{1}{2}} & X^*X\end{array}\right) \left(\begin{array}{clcr}P_1^{\frac{1}{2}}U_1P_1^{\frac{1}{2}} & P_1^{\frac{1}{2}}X\\0 & 0\end{array}\right) \\
  &=& \left(\begin{array}{clcr}P_1^{\frac{1}{2}}U^*_1P_1 & 0\\X^*P_1 & 0\end{array}\right) \left(\begin{array}{clcr}I_1 & U^*_1X\\X^*U_1 & X^*X\end{array}\right) \left(\begin{array}{clcr}P_1U_1P_1^{\frac{1}{2}} & P_1X\\0 & 0\end{array}\right) \\ 
  &=& \left(\begin{array}{clcr}P_1^{\frac{1}{2}}U^*_1P_1 & 0\\X^*P_1 & 0\end{array}\right) \left(\begin{array}{clcr}I_1 & 0\\X^*U_1 & I_2\end{array}\right) \left(\begin{array}{clcr}I_1 & U_1^*X\\0 & X^*U_1U^*_1X+X^*X\end{array}\right)\left(\begin{array}{clcr}P_1U_1P_1^{\frac{1}{2}} & P_1X\\0 & 0\end{array}\right)\\ 
  &=&   \left(\begin{array}{clcr}P_1^{\frac{1}{2}}U^*_1P_1\\X^*P_1 & 0\end{array}\right)\left(\begin{array}{clcr}P_1U_1P_1^{\frac{1}{2}} & P_1X\\0 & 0\end{array}\right)\\ 
  &=& A^*(P_1\oplus I_2)A=A^*Q^2_1A
\end{eqnarray*}
and
\begin{eqnarray*} B^*DB &=&  \left(\begin{array}{clcr}U^*_1P_1 & 0\\X^*P_1 & 0\end{array}\right) \left(\begin{array}{clcr}I_1 & U^*_1X\\X^*U_1 & X^*X\end{array}\right) \left(\begin{array}{clcr}P_1U_1 & P_1X\\0 & 0\end{array}\right) \\ 
&=& \left(\begin{array}{clcr}U^*_1P_1 & 0\\X^*P_1 & 0\end{array}\right) \left(\begin{array}{clcr}I_1 & 0\\X^*U_1 & I_2\end{array}\right) \left(\begin{array}{clcr}I_1 & U_1^*X\\0 & X^*U_1U^*_1X+X^*X\end{array}\right)\left(\begin{array}{clcr}P_1U_1 & P_1X\\0 & 0\end{array}\right)\\ 
&=& \left(\begin{array}{clcr}U^*_1P_1 & 0\\X^*P_1 & 0\end{array}\right)\left(\begin{array}{clcr}P_1U_1 & P_1X\\0 & 0\end{array}\right)\\ 
&=& B^*B. 
\end{eqnarray*}
Thus $T\in (m,|T^n|^2)$-expansive implies 
$$
\sum_{j=1}^m{(-1)^{j}\left(\begin{array}{clcr}m\\j\end{array}\right)A^{*j}Q_1^2A^j }+ C\leq 0$$ and $$\sum_{j=1}^m{(-1)^{j}\left(\begin{array}{clcr}m\\j\end{array}\right) B^{*j}B^j} +D\leq 0.
$$
Assume now that $\left(\begin{array}{clcr}I_1 & X\\X^* & X^*X\end{array}\right) \geq I$. Then
\begin{eqnarray*} 
D=\left(\begin{array}{clcr}I_1 &U_1^* X\\X^*U_1 & X^*X\end{array}\right) &=& (U^*_1\oplus I_2)\left(\begin{array}{clcr}I_1 & X\\X^* & X^*X\end{array}\right)(U_1\oplus I_2)\\ &\geq& (U^*_1\oplus I_2)(U_1\oplus I_2)=I
\end{eqnarray*} 
and

\begin{eqnarray*} 
 C=Q_1DQ_1=(P_1^{\frac{1}{2}}\oplus I_2)D(P_1^{\frac{1}{2}}+I_2)\geq (P_1+I_2).
 \end{eqnarray*}
Let $P_1\oplus I_2=Q^2_1=Q$; then $Q\geq 0$ is invertible and 
$$
||x||^2_Q=\langle x,x\rangle_Q=\langle Qx,x\rangle=||Q^{\frac{1}{2}}x||^2
$$ 
defines an equivalent norm on $\H$. Since
\begin{eqnarray*} 
T\in (m,|T^n|^2)-{\rm expansive} &\Longrightarrow&     \sum_{j=0}^m{(-1)^{j}\left(\begin{array}{clcr}m\\j\end{array}\right) B^{*j}B^j} +(D-I)\leq 0\\ &\Longleftrightarrow& \sum_{j=0}^m{(-1)^{j}\left(\begin{array}{clcr}m\\j\end{array}\right) B^{*j}B^j} \leq I-D\leq 0, 
\end{eqnarray*}
$B$ is $m$-expansive. Again, since
\begin{eqnarray*} 
T\in (m,|T^n|^2)-{\rm expansive} &\Longrightarrow&     \sum_{j=0}^m{(-1)^{j}\left(\begin{array}{clcr}m\\j\end{array}\right) A^{*j}QA^j} +(C-Q)\leq 0\\ &\Longleftrightarrow& \sum_{j=0}^m{(-1)^{j}\left(\begin{array}{clcr}m\\j\end{array}\right) A^{*j}QA^j} \leq Q-C\leq 0, 
\end{eqnarray*}
$A$ is $(m,Q)$-expansive. Equivalently, $A$ is $m$-expansive in $(\H, ||.||_Q)$.
\end{demo}

The case in which $T\in (2,P)$-expansive, and $P=T^*T$ with $0\notin\sigma(T^*T)$, is of some interest. Letting as before $T=\left(\begin{array}{clcr}T_1 & X\\0 & 0\end{array}\right)\in B({\overline{T(\H)}}\oplus T^{*-1}(0))$, it is seen that
$T_1$ is an invertible $(2,|T_1|^2)$-expansive operator. Consequently, $T_1^{*2}T_1^2=T_1^*T_1$ (see Remark \ref{rema01}), equivalently, $T_1$ is an invertible isometry (i.e., a unitary), and $\sigma(T)\subseteq \partial{\D}$ (where $\partial{\D}$ denotes the unit circle in the complex plane.) We note here
that $\sigma(T)\subseteq \sigma(T_1)\cup\{0\}$ and $0\notin\sigma_a(T)$ by hypothesis. This implies indeed that $T$ is a unitary, hence satisfies  almost every variety of Browder, Weyl theorems ({\em cf.} \cite[Lemma 3.13 and Theorem 3.14]{JKKL}); see \cite{A}, Chapters 5 and 6 for information on Browder, Weyl type theorems.


\vskip10pt \noindent\normalsize\noindent{B.P. Duggal\\ University of Ni\v s\\
Faculty of Sciences and Mathematics\\
P.O. Box 224, 18000 Ni\v s, Serbia}

\noindent {\it E-mail}: {\tt bpduggal@yahoo.co.uk}

\vskip6pt\noindent \noindent\normalsize\rm I. H. Kim\\Department of
Mathematics\\ Incheon National University\\ Incheon, 22012, Korea.\\
\noindent\normalsize \tt e-mail: ihkim@inu.ac.kr

\end{document}